\documentclass[12pt]{article}

\usepackage[T2A]{fontenc}
\usepackage[cp866]{inputenc}

\usepackage{amsfonts}
\usepackage{amssymb}
\usepackage{amsmath}
\usepackage{amsthm}

\def \le {\leqslant}
\def \ge {\geqslant}

\begin{document}

\centerline{\bf Density modulo 1 of sublacunary sequences: } \centerline{\bf application of
 Peres-Schlag's arguments}
\vskip+1.5cm \centerline{\bf Moshchevitin N.G. \footnote{ Research is supported by grants RFFI 06-01-00518, MD-3003.2006.1, NSh-1312.2006.1 and
INTAS 03-51-5070
 }
} \vskip+1.5cm

{\bf Abstract.}\,\, Let the sequence
 $\{t_n\}_{n=1}^{\infty}$ of reals satisfy the condition
$ \frac{t_{n+1}}{t_n} \ge 1+ \frac{\gamma}{n^\beta},\,\,\,0\le \beta <1,\,\,\, \gamma >0. $ Then   the set $ \{\,\, \alpha \in [0,1]:\,\,\,
\exists \varkappa > 0 \,\,\,\forall n \in \mathbb{N}
 \,\,\, ||t_n \alpha ||
> \frac{\varkappa }{n^\beta \log (n+1)} \,\,
\} $ is uncountable. Moreover its Hausdorff dimension is equal to
1.
 Consider the set of naturals of the form $2^n3^m$ and let the sequence
$ s_1{=}1, s_2{=}2, s_3{=}3, s_4{=}4, s_5{=}6, s_6 = 8,\dots $ performs this set as an increasing sequence.  Then the set
$
\{\,\, \alpha \in [0,1]:\,\,\, \exists \varkappa > 0 \,\,\,\forall n \in \mathbb{N}
 \,\,\, ||s_n \alpha ||
> \frac{\varkappa}{\sqrt{n}\log (n+1)} \,\,
\} $ also has Hausdorff dimension equal to  1. The results
obtained use an original approach due to Y. Peres and W. Schlag.

 \vskip+1.5cm

{\bf 1. Introduction.} \,\,\, A sequence $\{t_j\},\,\,\ j = 1,2,3,..$ of positive real numbers is defined to be lacunary if for some $ M>0$ one
has
$$
\frac{t_{j+1}}{t_j} \ge 1+\frac{1}{M},\,\,\, \forall j \in \mathbb{N}.
$$
 Erd\"os \cite{E} conjectured that for any lacunary sequence there exists
real $\alpha$ such that the set of fractional parts $ \{ \alpha t_j\} , \,\, j \in\mathbb{N}$ is not dense in $[0,1]$. This conjecture was
proved by A. Pollington \cite{Po} and B. de Mathan \cite{M}. Some quantitative improvements were due to Y. Katznelson \cite{K}, R. Akhunzhanov
and N. Moshchevitin \cite{AM} and A. Dubickas \cite{D}. The best known quantitative estimate is due to Y. Peres and W. Schlag \cite{P}. The last
authors proved that with some positive constant $\gamma > 0$ for any sequence $\{t_j\}$ under  consideration there exists a real number $\alpha
$ such that
$$
||\alpha t_j  || \ge \frac{\gamma}{M\log M},\,\,\, \forall j \in \mathbb{N}.
$$
Y. Peres and W. Schlag use an original approach connected with the Lovasz local lemma.

From another hand R. Akhunzhanov and N. Moshchevitin in \cite{MAH} generalized Pollington - de Mathan's result to sublacunary sequences. For
example for a sequence  $\{t_j\}$ under condition
$$
\frac{t_{j+1}}{t_j} \ge 1+\frac{\gamma}{n^\beta},\,\,\, \forall j \in \mathbb{N}, \gamma > 0, \beta \in (0,1/2]
$$
they proved the existence of real irrational $\alpha$ such that
$$
\liminf_{n\to \infty}\left( ||t_n \alpha || \times n^{2\beta }\right) > 0.
$$
Another application from \cite{MAH} deals with the sequence of
naturals of the form $ 2^m3^n,\,\, m,n \in\mathbb{N}\cup \{ 0\}$.

In the present paper we apply the arguments from \cite{P} to improve the results from \cite{MAH} mentioned above.

 {\bf 2. Results.}\,\,\,
Let $1 \le t_1< t_2 < ... t_n < t_{n+1} < ...$ be a strictly increasing sequence of reals and $ \lim_{n\to \infty} t_n = +\infty$. For a given
sequence $\{t_n\}$ we define the function
\begin{equation}
H(n,\tau ) = \min \left\{ k\in \mathbb{N}:\,\,\, \frac{t_{n+k}}{t_n} \ge \tau \right\}. \label{H}
\end{equation}

 {\bf Theorem 1.}\,\,\,{\it Let $ 0< \eta < 1$.
Consider a sequence $\{ h(n)\}_{n=1}^\infty \subset \mathbb{N}$ of natural numbers such that
 for all natural $n$ under condition $ n
> h(n)$
 the function $ n\mapsto n - h(n)$ is increasing and
a decreasing sequence $\{\delta (n) \}_{n=1}^{\infty} $ of positive real numbers. Let the sequence $\{n_k\}_{k=0}^K$ of natural numbers is
defined to satisfy the condition
\begin{equation}
n_k = n_{k+1}-h(n_{k+1}) \label{K}
\end{equation}
for $  0\le k \le K-1$.
 Let
our sequences
  satisfy
the following conditions (i),(ii) and (iii) below.

(i)  For any natural $n$ under condition $ n > h(n)$ the following inequality is valid
$$
h(n) \ge H(n-h(n), 1/\delta (n-h(n))).
$$

(ii)  For any   $k\le K -1$ the following inequality is valid
$$
\sum_{v = n_k+1}^{n_{k+1}-1} \delta ( v ) \le  \frac{(1-\eta)\eta}{4}.
$$

(iii) For $k=0$ the following inequality is valid
$$
\sum_{v = 1}^{n_{0}} \delta ( v ) \le  \frac{1-\eta}{16}.
$$
 Then
for the set
$$
{\cal A}_K  = \{ \alpha \in [0,1]:\,\,\, ||t_n \alpha ||
> \delta (n) \,\,
\forall n \le n_K \}
$$
one has
$$
\mu ({\cal A}_K ) \ge \eta^{K+1}.
$$}
Here $\mu (\cdot ) $ denotes the  Lebesgue measure. Note that the
sets $ {\cal A}_K $ are closed and nested: $ {\cal A}_{K+1}
\subseteq {\cal A}_K.$ Moreover if we have a natural number $N$ we
can construct a sequence  $ \{n_k\}$ such that $ n_K = N$, the
equalities (\ref{K}) are satisfied,  $ n_0 = n_1 - h(n_1) \ge 1$
but $ n_0 - h(n_0) \le 0$. Hence as a corollary of Theorem 1 we
immediately obtain

 {\bf Theorem 2.}\,\,\,{\it Let $ 0< \eta < 1$.
Consider a sequence $\{ h(n)\}_{n=1}^\infty \subset \mathbb{N}$ of
natural numbers such that for all natural $n$ under condition $ n
> h(n)$
the function $ n\mapsto n - h(n)$ is increasing and a decreasing sequence $\{\delta (n) \}_{n=1}^{\infty} $ of positive reals.
 Let
these sequences
  satisfy
the following conditions (i) from Theorem 1 and the conditions (ii
$'$) (iii$'$) below.

(ii $'$)  For all natural numbers $n$ under condition  $ n > h(n) $ the following inequality is valid
$$
\sum_{v = n - h(n) +1}^{n-1} \delta ( v ) \le  \frac{(1-\eta)\eta}{4}.
$$

(iii $'$)   For all natural numbers $n$ under condition  $ n \le h(n) $ the following inequality is valid
$$
\sum_{v = 1}^{n} \delta ( v ) \le  \frac{1-\eta}{16}.
$$
Then the set
$$
{\cal A}  = \{ \alpha \in [0,1]:\,\,\, ||t_n \alpha ||
> \delta (n) \,\,
\forall n \in \mathbb{N}\}
$$
is nonempty. }

{\bf Theorem 3.}\,\,\,{\it Let the conditions of theorem 2 be
satisfied and an  infinite sequence $\{n_k\}_{k=0}^\infty$ of
naturals  satisfies the condition (\ref{K}) for all natural $k$.
 Let the series
\begin{equation}
\sum_{k=1}^\infty \frac{1}{ \eta^k} \cdot \left(
\frac{t_{n_k}}{\delta (n_k)}\right)^\nu /
\left(\frac{t_{n_{k-1}}}{\delta (n_{k-1})}\right)
 \label{ser} \end{equation}
converges for all $ \nu < \nu_0$ Then the set $ {\cal A} $ from
Theorem 2 has Hausdorff dimension   $ \ge \nu_0$. }

We give a complete proof of theorem 1 in Sections 3,4. In Section 5 we give comments to the proof of Theorem 3. In section 6 we give some
applications of our results.

{\bf 4. Lemmata.}\,\,\, For $ n \ge 1 $ we define
\begin{equation}
l_n = \left\lfloor \log_2 \left(\frac{t_n}{2\delta (n)}\right) \right\rfloor. \label{l}
\end{equation}
From monotonicity of $t_n$ and $\delta (n)$ it follows that $ l_{n+1} \ge l_n$. Put
$$
E(n,a) = \left[ \frac{a}{t_n} - \frac{\delta (n)}{t_n}, \frac{a}{t_n} + \frac{\delta (n)}{t_n} \right]
$$
Let $A_n$ be the union of dyadic intervals of the form
$$
\left( \frac{b}{2^{l_n}},\frac{b+\varepsilon}{2^{l_n} } \right),\,\,\, b \in \mathbb{Z},\,\, \varepsilon \in \{ 1,2\}
$$
which covers the set
$$
\bigcup_{0\le a \le \lceil t_n\rceil} E (n,a) \bigcap [0,1].
$$
So
$$
\bigcup_{0\le a \le \lceil t_n\rceil} E (n,a) \bigcap [0,1]
\subseteq A_n. $$ Define $A_n^c = [0,1]\setminus A_n$. Note that
$$
\mu ( A_n) \le (\lceil t_n \rceil +1) \frac{2\delta (n)}{t_n} \le 16\delta (n)
$$
and
\begin{equation}
\mu \left( \bigcap_{n\le n_0} A_n^c\right) \ge 1 - 16 \sum_{n=1}^{n_0} \delta (n). \label{init}
\end{equation}

 {\bf Lemma 1.}\,\,\,{\it Let $ n > h(n)$.
 Let the condition (i) holds and
$$
\mu \left(\bigcap_{j \le n- h(n)} A_j^c\right) > 0.
$$
Then \begin{equation} \mu \left(\bigcap_{j \le n- h(n)} A_j^c\bigcap A_{n}\right) \le 4\delta (n)
 \mu \left(\bigcap_{j \le n-h(n)} A_j^c\right).
\label{yyy} \end{equation} }

Proof. The set $\bigcap_{j \le n- h(n)} A_j^c$ can be considered as  a union \begin{equation} \bigcap_{j \le n-h(n)} A_j^c = \bigcup_{\nu = 1}^T
I_\nu \label{I}
\end{equation}
 of the dyadic  intervals $I_\nu = I_\nu ^{(n-h(n))}$ of the form
$$
\left[ \frac{b}{2^{l_{n-h(n)}}}, \frac{b+1}{2^{l_{n-h(n)}}} \right],\,\,\, b \in \mathbb{Z}
$$
where $T\ge 1$. Now the set $ A_{n}\cap I_\nu$ can be represented as a union
$$
A_{n}\bigcap I_\nu = \bigcup_{i = 1}^{W_\nu} J_i
$$
of intervals $J_i$ of the form
$$
\left[ \frac{b}{2^{l_{n)}}}, \frac{b+1}{2^{l_{n}}} \right] .
$$
Moreover
$$
W_\nu \le\left\lfloor \left( \frac{1}{2^{l_{n-h(n)}}} + \frac{\delta (n)}{2^{l_{n}}}\right) t_{n}\right\rfloor+1 \le
\frac{t_{n}}{2^{l_{n-h(n)}}} + 2 .
$$
So
$$
\mu \left( A_{n}\cap I_\nu \right) = \frac{W_\nu}{2^{l_{n}}}
$$
and
$$
\mu \left(\bigcap_{j \le n-h(n)} A_j^c\bigcap A_{n}\right) \le \frac{T}{2^{l_{n}}}\left( \frac{t_{n}}{2^{l_{n-h(n)}}}+2\right)= \mu
\left(\bigcap_{j \le n-h(n)} A_j^c\right)\frac{2^{l_{n-h(n)}}}{2^{l_{n}}}\left( \frac{t_{n}}{2^{l_{n-h(n)}}}+2\right)=
$$
$$
= \mu \left(\bigcap_{j \le n-h(n)} A_j^c\right) \left( \frac{t_{n}}{2^{l_{n}}} + 2\cdot \frac{2^{l_{n-h(n)}}}{2^{l_{n}}} \right).
$$
But
\begin{equation}
\frac{t_{n}}{2^{l_{n}}}\le 2 \delta (n ) \label{1}
\end{equation}
 from the definition of $l_n$ (formula (\ref{l})). For the second summand we
have
\begin{equation}
\frac{2^{l_{n-h(n)}}}{2^{l_{n}}} \le 2\cdot
\frac{t_{n-h(n)}}{t_{n}}\cdot  \frac{\delta (n)}{\delta (n-h(n))}
 \le 2 \delta (n) \label{2}
\end{equation}
 from the condition (i) and the definition
(\ref{H}) of the function $H(\cdot , \cdot )$.

Now Lemma 1 follows from (\ref{1},\ref{2}).

For fixed $\tau$ and $ 0\le v \le h(\tau)$ define $\tau_v = \tau -h(\tau) +v$. Note that $ \tau_{h(\tau)} = \tau$ and $ \tau_0 = \tau -
h(\tau)$. Note that
$
\tau_0 \le \tau_v \le\tau .$

 {\bf Lemma 2.}\,\,\, {\it Let   the function $ n-h(n)$ is increasing
and the condition (i) holds. Let for $ \tau_0>h (\tau_0 ) $ the
following inequality is valid:
\begin{equation}
\mu\left( \bigcap_{j\le \tau_0} A_j^c\right) \ge \eta \mu\left( \bigcap_{j\le \tau_0- h(\tau_0)} A_j^c\right)
>0
\label{iii}
\end{equation}
with some positive $\eta$.

 Then
 we have
 \begin{equation}
\mu\left( \bigcap_{j\le \tau} A_j^c\right) \ge \left( 1 - \frac{4}{\eta }\,\, \sum_{v=\tau_1}^{\tau-1}\delta (v) \right) \times \mu\left(
\bigcap_{j\le \tau_0} A_j^c\right) . \label{l2} \end{equation} }

Proof.

We have
$$
\mu\left( \bigcap_{j\le \tau} A_j^c\right) = \mu\left( \left(\cdots \left(\left( \bigcap_{j \le \tau- h(\tau)} A_j^c\right) \setminus
A_{\tau-h(\tau) +1}\right) \setminus \cdots \right) \setminus A_{\tau}\right) \ge
$$
$$
\ge \mu\left( \bigcap_{j\le \tau-h(\tau)} A_j^c\right) - \sum_{v =1}^{h(\tau) } \mu \left( A_{\tau_v} \bigcap\left(\bigcap_{j\le \tau - h(\tau)}
A_j^c\right)\right).
$$
But as $\tau_v \le \tau$ from the monotonicity condition for $ n - h(n)$   we get $ \tau - h(\tau) \ge \tau_v - h(\tau_v)$ so
\begin{equation}
\bigcap_{j\le \tau - h(\tau)} A_j^c\subseteq \bigcap_{j\le \tau_v - h(\tau_v)} A_j^c . \label{mono} \end{equation}
 Now
$$
\mu\left( \bigcap_{j\le \tau} A_j^c\right) \ge \mu\left( \bigcap_{j\le \tau - h(\tau)} A_j^c\right) - \sum_{v =1}^{h(\tau) } \mu \left(
A_{\tau_v} \bigcap\left(\bigcap_{j\le \tau_v - h(\tau_v)} A_j^c\right)\right).
$$
We apply Lemma 1 for $ n = \tau_v,\,\, v = 1,...,h(\tau) $ (it is possible as from (\ref{mono}) and $ \mu\left( \bigcap_{j\le \tau- h(\tau)}
A_j^c\right)
>0
$ it follows that $ \mu\left( \bigcap_{j\le \tau_v- h(\tau_v)} A_j^c\right)
>0
$ for all $v$)
 and obtain the inequality
$$
\mu \left(A_{\tau_v} \bigcap\left(\bigcap_{j\le \tau_v - h(\tau_v)} A_j^c\right)\right) \le 4 \delta (\tau_v) \mu \left( \bigcap_{j\le \tau_v -
h(\tau_v)} A_j^c\right) .
$$
Now
$$
\mu\left( \bigcap_{j\le \tau} A_j^c\right) \ge \mu\left( \bigcap_{j\le \tau - h(\tau)} A_j^c\right) - 4\left(\sum_{v =1}^{h(\tau) } \delta
(\tau_v) \right)
 \times
\max_{1\le v < h(\tau) } \mu \left( \bigcap_{j\le \tau_v - h(\tau_v)} A_j^c\right) \ge
$$
$$
\mu\left( \bigcap_{j\le \tau - h(\tau)} A_j^c\right) - 4\left(\sum_{v =1}^{h(\tau) } \delta (\tau_v) \right)
 \times
\max_{0\le v < h(\tau) } \mu \left( \bigcap_{j\le \tau_v - h(\tau_v)} A_j^c\right) .$$
 But  we have the condition that the function $ n-h(n)$ is
increasing. So the maximum here is obtained at $ v = 0$. It follows that
$$
\mu\left( \bigcap_{j\le \tau} A_j^c\right) \ge \mu\left( \bigcap_{j\le \tau - h(\tau)} A_j^c\right) - 4\left( \sum_{v =1}^{h(\tau) } \delta
(\tau_v) \right) \times \mu  \left(\bigcap_{j\le \tau_0 - h(\tau_0)} A_j^c\right) .
$$
We apply (\ref{iii}) below:
$$
\mu\left( \bigcap_{j\le \tau} A_j^c\right) \ge \mu\left( \bigcap_{j\le \tau - h(\tau)} A_j^c\right) - \frac{4}{\eta}\,\, \left( \sum_{v
=1}^{h(\tau) } \delta (\tau_v) \right) \times\mu  \left(\bigcap_{j\le \tau_0 } A_j^c\right) .
$$
 Remember that $\tau_0 = \tau-h(\tau)$ and
 Lemma 2 follows.

{\bf 4. Proof of Theorem 1.} \,\,\, From condition (iii) of the Theorem 1 and (\ref{init}) it follows that $ \mu\left( \bigcap_{j\le n_0}
A_j^c\right) \ge \eta \ge \eta \mu\left( \bigcap_{j\le n_0-h(n_0)} A_j^c\right)$. This is the base of induction. The inductive step $\mu\left(
\bigcap_{j\le n_{k+1}} A_j^c\right) \ge \eta \mu\left( \bigcap_{j\le n_k} A_j^c\right)$ follows from condition (ii) and Lemma 2: We must put $
\tau = n_{k+1}$, then $ \tau_0 = n_k$. From inductive hypothesis we have (\ref{iii}). The condition (ii) leads to inequality $ 1 - \frac{4}{\eta
}\,\, \sum_{v=\tau_1}^{\tau-1}\delta (v) \ge \eta . $

{\bf 4. Sketched proof of  of Theorem 3.}\,\,\, In order to prove Theorem 3 one must do the following. In the proof of Theorem 1 instead of the
inequality (\ref{yyy}) of Lemma 1 one should prove
$$
\mu \left( I_\nu^{(n-h(n))} \bigcap A_{n}\right) \le 4\delta (n) \mu \left(I_\nu^{(n-h(n))}\right),
$$
where $I_\nu^{(n-h(n))} $ is from  partition (\ref{I}). Then under the condition
$$
\mu\left( I_{\nu '}^{( \tau_0- h(\tau_0))}\cap A_{\tau_0}\right) \ge \eta \mu\left( I_\nu^{( \tau_0- h(\tau_0))} \right)
>0
$$
 one should prove instead of the inequality (\ref{l2}) of Lemma 2 the following
inequality:
$$
\mu\left( I_{\nu }^{(\tau_0)}\bigcap \left( \bigcap_{j\le \tau }A_{j}^c \right) \right) \ge \left( 1 - \frac{4}{\eta }\,\,
\sum_{v=\tau_1}^{\tau-1}\delta (v) \right) \times \mu\left( I_{\nu } ^{(\tau_0)} \right).
$$
It means that in each interval of the form $ I^{(\tau_0)}_\nu$ there exist not less than
$$
N= \frac{\mu\left( I_{\nu }^{(\tau_0)}\bigcap \left( \bigcap_{j\le \tau }A_{j}^c \right) \right) }{\mu\left( I_{\nu '} ^{(\tau)} \right) }\ge
\eta 2^{l_\tau - l_{\tau_0}}
$$
pairwise  disjoint subintervals of the form $I_{\nu '} ^{(\tau)}$.
 Then as in \cite{MAH} one should take into account the convergence of (\ref{ser}) and apply the following well-known result:

 {\bf Theorem (Eggleston \cite{Eg}).}\,\,\ {\it Let for every $k$ we have a set
$A_k{=}\bigsqcup\limits_{i=1}^{R_k}I_{k}(i)$ where  $I_{k}(i)$ are segments of real line of length $|I_{k}(i)|= \Delta_k$. Let each interval
$I_{k}(i)$ has exactly $N_{k+1}{>}1$ pairwise disjoint subintervals $I_{k+1}(i')$  of length $\Delta_{k+1}$ from the set
 $A_{k+1}$. Let
$R_{k+1}{=}R_k{\cdot} N_{k+1}$ . Suppose
 $0{<}\nu_0{\le}1$and for every $0{<}\nu{<}\nu_0$  the series
$ \sum_{k=2}^{\infty}\frac{\Delta_{k-1}}{\Delta_k}(R_k(\Delta_k)^\nu)^{-1} $ converges. Then the set $A{=}\bigcap_{k=1}^{\infty}A_k$ has
Hausdorff dimension ${\rm HD}(A){\ge}\nu_0$. }

 {\bf 6.  Examples.}\,\,\,
Note that the proof of Theorem 1 follows directly the arguments by Y.Peres and W. Schlag from \cite{P}. The author in \cite{MPRE} (following
Peres-Schlag's arguments)
 established for lacunary
sequence $\{t_n\}$ under condition
$$
\frac{t_{j+1}}{t_j} \ge 1+\frac{1}{M},\,\,\, \forall j \in \mathbb{N}.
$$
 the existence
of a real number $ \alpha$ such that
$$
||\alpha t_j|| \ge \frac{1}{2^{11}M\log M},\,\,\, \forall j \in \mathbb{N}.
$$.
We consider  some examples with sublacunary sequences below.

{\bf A. \,\,\, Sublacunary sequences.} \,\,\, Let $\{t_n\}_{n=1}^{\infty}$ satisfy the condition
\begin{equation}
  \frac{t_{n+1}}{t_n}  \ge 1+
\frac{\gamma}{n^\beta},\,\,\, 0\le \beta <1,\,\,\, \gamma >0 .\label{doad}
\end{equation}
 We take $\eta < 1$ close to 1 and
\begin{equation}
  h(n) = \lfloor c_1n^\beta\log (n+c_2) \rfloor ,
 \,\,\,
\delta (n) =    \frac{(1-\beta)(1-\eta)\eta}{2^5 c_1 (n+c_2)^\beta
\log ( n+c_2)}, \label{para} \end{equation}
 Here large positive constants
 $c_1,c_2$  (depending on $\beta $ and $\eta $)
 should be defined in the following way.
 In our situation under condition $ n > h(n)$ for $ \gamma_1< \gamma_1$ one has
$$
\frac{t_{n}}{t_{n-h(n)}} \ge \prod_{j=n-h(n)}^{n-1} \left( 1+ \frac{\gamma}{j^\beta} \right) \ge \exp \left( \sum_{j=n-h(n)}^{n-1} \log \left(
1+\frac{\gamma }{j^\beta }\right)\right) \ge \exp \left(\omega \frac{h(n)}{n^\beta}\right)\ge (n+c_2)^{\omega c_1}
$$
with $ \omega = \omega (\beta, \gamma_1 )$.
 Let $c_1= c_1(\beta, \eta  )$ be a large positive constant
such that for all real $  y \ge 2$ we have
$$
y^{\omega c_1} \ge \frac{2^5 c_1 y^\beta \log y}{(1-\beta)(1-\eta)
\eta}.
$$
Then
$$
\frac{t_{n}}{t_{n-h(n)}} \ge (n+c_2)^{\omega c_1}\ge \frac{
 2^5 c_1 (n+c_2)^\beta \log (n+c_2)}{(1-\beta)(1-\eta)\eta} =\frac{1}{\delta(n)}\ge \frac{1}{\delta(n-h(n))}
 $$
and  the condition (i$'$) of Theorem 2 is satisfied.

 So we  have $c_1$ fixed and then we  define $c_2$.
 Let $c_2 = c_2(\beta )$ be a large positive constant such that
\begin{equation}
\max_{n \in \mathbb{N}} \frac{4c_1\log (n+c_2)}{(n+c_2)^{1-\beta }} \le 1, \label{PP}
\end{equation}
\begin{equation}
h\left(\frac{1}{2^5\delta (0)}\right) = \left\lfloor c_1\left(
\frac{c_1c_2^\beta\log c_2}{(1-\beta)(1-\eta )\eta}\right)^\beta
\log \left( \frac{2^{2}c_1c_2^\beta\log c_2}{1-\beta} + c_2\right)
 \right\rfloor\le
 \frac{1}{2^5\delta (0)} =
\frac{c_1c_2^\beta\log c_2}{(1-\beta)(1-\eta )\eta}.
 \label{QQ}
\end{equation}
\begin{equation}
\min_{y \ge 1}\left( (1-\beta ) \log ( y + c_2)-\frac{y}{y+c_2} \right) > 0\label{RRR}
\end{equation}
 Then from (\ref{PP}) it follows that $ \frac{h(n)}{n+c_2}\le \frac{1}{2}$ and for $ n >
h(n)$ we have
$$
\sum_{v = n - h(n) +1}^{n-1} \delta ( v ) \le
\frac{(1-\beta)(1-\eta)\eta}{2^5 c_1 \log ( n-h(n)+c_2)} \sum _{v
= n - h(n) +1}^{n-1} \frac{1}{v^\beta} \le (1-\eta)\eta \times
\frac{n^{1-\beta} - (n-h(n))^{1-\beta}}{ 2^4 c_1 \log (
n-h(n)+c_2)}\le
$$
$$
\le \frac{(1-\eta)\eta h(n)}{2^4 c_1 n^\beta \log ( n-h(n)+c_2)}
\le   \frac{ (1-\eta)\eta\log (n+c_2)}{ 2^3  \log ( n-h(n)+c_2)}=
$$
$$=
\frac{ (1-\eta)\eta }{ 2^3 } \times \frac{ \log (n+c_2)}{ \log (
n+c_2) + \log ( 1- \frac{h(n)}{n+c_2})}\le \frac{ (1-\eta)\eta }{
2^3 } \times \frac{ \log (n+c_2)}{ \log ( n+c_2) -\log 2} \le
\frac{(1-\eta)\eta}{4}.
$$
So the condition (ii$'$) of Theorem 2 is satisfied.

Moreover for the value $ n_0 = n_0(\beta, c_1,c_2 ) = \max\{ n \in \mathbb{N}:\,\,\, n \le h (n) \} $ from (\ref{QQ}) it follows that $n_0 \le
\frac{1}{2^5\delta (0)}$ and
 the condition (iii$'$) of Theorem 2 is satisfied also.

 Also we must note that if $ y \ge 1$ and $ y > h(y) \ge
c_1y^\beta \log ( y + c_2)$ then the function $ y - c_1 y^\beta
\log ( y + c_2)$ is increasing as from (\ref{RRR}) it follows that
$$
(y - c_1 y^\beta \log ( y + c_2))' = 1-\beta c_1y^{\beta - 1} \log (y+c_2) - \frac{c_1y^\beta}{y+c_2} = c_1 y^{\beta - 1}\left( (1-\beta ) \log
( y + c_2)-\frac{y}{y+c_2}\right)> 0.
$$

Now we have checked all the conditions of Theorem 2. It follows
that the set
$${\cal B}=  \{\,\, \alpha \in [0,1]:\,\,\, \exists \varkappa > 0 \,\,\,\forall n \in
\mathbb{N}
 \,\,\, ||t_n \alpha ||
> \frac{\varkappa }{n^\beta \log (n+1)} \,\,
\} $$ is nonempty (obviously, uncountable).

 Note that
the set $\{ n\in \mathbb{N}:\,\,\, n\le h(n)\}$ is finite. Hence  we can construct a sequence of naturals  $\{n _k\}$ satisfying (\ref{K}).

If it happens that in addition to (\ref{doad}) we have
\begin{equation}
  \frac{t_{n+1}}{t_n} \le 1+
\frac{\gamma _2}{n^\beta} \label{doadd}
\end{equation}
with
 some $ \gamma_2 >\gamma$
then for the sequence $\{n _k\}$ we get $t_{n_k} \le t_{n_{k-1}}
n^{\gamma_3} $ and $ k\le \gamma_4 n_k^{1-\beta}$ with positive
$\gamma_{3,4}$. Now
$$
\frac{1}{\eta^k} \cdot \frac{t_{n_k}^\nu}{t_{n_{k-1}}} \ll
\frac{1}{e^{\gamma_5 n_k^{1-\beta}}}\cdot \frac{1}{\eta^k} \ll
\frac{1}{(e^{\gamma_5}\eta^{\gamma_5})^{n^{1-\beta}}}
$$
(here all constants  $\gamma_j$ do not depend on $\eta$) and for
$\eta$ close to 1 the series (\ref{ser}) converges. From Theorem 3
it follows that
 the set
$ {\cal B}$
 has Hausdorff dimension equal to  1.
We should note that it is possible to choose function $h(n)$
(actually in the same manner as it was done in \cite{MAH}) to
satisfy the conditions of Theorem 3 without additional assumption
(\ref{doadd}) on the rate of growth of the sequence $t_n$.

We should note that it would be interesting to investigate {\it winning} properties of the considered sets (for the definition of winning sets
see \cite{S},\cite{SS}, for some partial results see \cite{MW}).

{\bf B. \,\,\, Subexponentional sequences.} \,\,\, Let $\{t_n\}_{n=1}^{\infty}$ satisfy the condition
\begin{equation}
\gamma_1\exp(  n^\beta) \le
  {t_n}
\le \gamma_2\exp(  n^\beta) ,\,\,\, 0 < \beta <  1,\,\,\, \gamma_{1,2} >0 .\label{doadexp}
\end{equation}
Then by the same reasons (as in example {\bf A}) we have that the Hausdorff dimension of the set
$$  \{\,\, \alpha \in [0,1]:\,\,\, \exists \varkappa > 0 \,\,\,\forall n \in
\mathbb{N}
 \,\,\, ||t_n \alpha ||
> \frac{\varkappa }{n^{1-\beta} \log (n+1)} \,\,
\} $$ is equal to  1.

 {\bf C. \,\,\, F\"urstenberg's sequence.}\,\,\, Consider the set of naturals of the
form $2^n3^m$ and let the sequence
$$
s_1{=}1, s_2{=}2, s_3{=}3, s_4{=}4, s_5{=}6, s_6 = 8,\dots
$$
performs this set as an increasing sequence. F\"urstenberg  \cite{FU} (see also \cite{BOSH}) proved that for any irrational
 $\alpha$ the set
of fractional parts $\{ 2^n3^m\alpha\}$ is dense in $[0,1]$. Hence
$$
\liminf_{n\to\infty}||s_n \alpha||{=}0.
$$
We should note that we no nothing about the rate of convergence to zero here. Obviously for $\alpha = 1/5$ one has
$$
||s_n/5|| \ge 1/5.$$ But $1/5$ is a rational number.

The sequence $\{ s_n\}$  satisfy (\ref{doadexp}) with $\beta = 1/2$.
  So from example {\bf B} it follows that Hausdorff dimension of the set
$$
\{\,\, \alpha \in [0,1]:\,\,\, \exists \varkappa > 0 \,\,\,\forall n \in \mathbb{N}
 \,\,\, ||s_n \alpha ||
> \frac{\varkappa}{\sqrt{n}\log (n+1)} \,\,
\}
$$
is equal to  1.

\vskip+2.0cm

author: Nikolay Moshchevitin

\vskip+0.5cm

e-mail: moshchevitin@mech.math.msu.su, moshchevitin@rambler.ru

\end{document}